\documentclass[11pt]{amsart} 
\usepackage{latexsym, amsmath, amsthm, amssymb, setspace, verbatim}
\usepackage[all]{xy}
\usepackage[applemac]{inputenc}
\usepackage{hyperref}

\title[The Rokhlin property vs. Rokhlin dimension 1]{ The Rokhlin property vs. Rokhlin dimension 1 on unital Kirchberg algebras }
\author[S. Barlak, D. Enders, H. Matui, G. Szabó, W. Winter]{ Sel\c{c}uk Barlak, Dominic Enders, Hiroki Matui, Gábor Szabó \and Wilhelm Winter }

\address{Department of Mathematics and Computer Science, University of Southern \linebreak \text{}\hspace{3mm} Denmark, Campusvej 55, 5230 Odense M, Denmark}
\email{barlak@imada.sdu.dk}

\address{Department of Mathematical Sciences,
Universitetsparken 5, \phantom{----------------}\linebreak \text{}\hspace{3.5mm}
2100 K{\o}penhavn,  Denmark}
\email{dominic.enders@math.ku.dk}

\address{Graduate School of Science, Chiba University,\phantom{------------------------------------}\linebreak \text{}
\hspace{2.3mm} Inage-ku, Chiba 263-8522, Japan}
\email{matui@math.s.chiba-u.ac.jp}

\address{Westfälische Wilhelms-Universität, Fachbereich Mathematik, \phantom{---------------}\linebreak \text{}\hspace{3.5mm} Einsteinstrasse 62, 48149 Münster, Germany}
\email{gabor.szabo@uni-muenster.de, wwinter@uni-muenster.de}

\thanks{\emph{Supported by:} SFB 878 \emph{Groups, Geometry and Actions}, GIF Grant 1137-30.6/2011 and \\
\emph{Danish National Research Foundation through the Centre for Symmetry and Deformation} (DNRF92)}

\begin{document}

% math
\renewcommand\matrix[1]{\left(\begin{array}{*{10}{c}} #1 \end{array}\right)}  % Matrix
\newcommand\set[1]{\left\{#1\right\}}  % Menge
\newcommand\mset[1]{\left\{\!\!\left\{#1\right\}\!\!\right\}}

%% Besondere Variablen
%Zahlmengen-Stil
\newcommand{\IA}[0]{\mathbb{A}} \newcommand{\IB}[0]{\mathbb{B}}
\newcommand{\IC}[0]{\mathbb{C}} \newcommand{\ID}[0]{\mathbb{D}}
\newcommand{\IE}[0]{\mathbb{E}} \newcommand{\IF}[0]{\mathbb{F}}
\newcommand{\IG}[0]{\mathbb{G}} \newcommand{\IH}[0]{\mathbb{H}}
\newcommand{\II}[0]{\mathbb{I}} \renewcommand{\IJ}[0]{\mathbb{J}}
\newcommand{\IK}[0]{\mathbb{K}} \newcommand{\IL}[0]{\mathbb{L}}
\newcommand{\IM}[0]{\mathbb{M}} \newcommand{\IN}[0]{\mathbb{N}}
\newcommand{\IO}[0]{\mathbb{O}} \newcommand{\IP}[0]{\mathbb{P}}
\newcommand{\IQ}[0]{\mathbb{Q}} \newcommand{\IR}[0]{\mathbb{R}}
\newcommand{\IS}[0]{\mathbb{S}} \newcommand{\IT}[0]{\mathbb{T}}
\newcommand{\IU}[0]{\mathbb{U}} \newcommand{\IV}[0]{\mathbb{V}}
\newcommand{\IW}[0]{\mathbb{W}} \newcommand{\IX}[0]{\mathbb{X}}
\newcommand{\IY}[0]{\mathbb{Y}} \newcommand{\IZ}[0]{\mathbb{Z}}

%Geschwungener Stil
\newcommand{\CA}[0]{\mathcal{A}} \newcommand{\CB}[0]{\mathcal{B}}
\newcommand{\CC}[0]{\mathcal{C}} \newcommand{\CD}[0]{\mathcal{D}}
\newcommand{\CE}[0]{\mathcal{E}} \newcommand{\CF}[0]{\mathcal{F}}
\newcommand{\CG}[0]{\mathcal{G}} \newcommand{\CH}[0]{\mathcal{H}}
\newcommand{\CI}[0]{\mathcal{I}} \newcommand{\CJ}[0]{\mathcal{J}}
\newcommand{\CK}[0]{\mathcal{K}} \newcommand{\CL}[0]{\mathcal{L}}
\newcommand{\CM}[0]{\mathcal{M}} \newcommand{\CN}[0]{\mathcal{N}}
\newcommand{\CO}[0]{\mathcal{O}} \newcommand{\CP}[0]{\mathcal{P}}
\newcommand{\CQ}[0]{\mathcal{Q}} \newcommand{\CR}[0]{\mathcal{R}}
\newcommand{\CS}[0]{\mathcal{S}} \newcommand{\CT}[0]{\mathcal{T}}
\newcommand{\CU}[0]{\mathcal{U}} \newcommand{\CV}[0]{\mathcal{V}}
\newcommand{\CW}[0]{\mathcal{W}} \newcommand{\CX}[0]{\mathcal{X}}
\newcommand{\CY}[0]{\mathcal{Y}} \newcommand{\CZ}[0]{\mathcal{Z}}

%Script Stil
\newcommand{\FA}[0]{\mathfrak{A}} \newcommand{\FB}[0]{\mathfrak{B}}
\newcommand{\FC}[0]{\mathfrak{C}} \newcommand{\FD}[0]{\mathfrak{D}}
\newcommand{\FE}[0]{\mathfrak{E}} \newcommand{\FF}[0]{\mathfrak{F}}
\newcommand{\FG}[0]{\mathfrak{G}} \newcommand{\FH}[0]{\mathfrak{H}}
\newcommand{\FI}[0]{\mathfrak{I}} \newcommand{\FJ}[0]{\mathfrak{J}}
\newcommand{\FK}[0]{\mathfrak{K}} \newcommand{\FL}[0]{\mathfrak{L}}
\newcommand{\FM}[0]{\mathfrak{M}} \newcommand{\FN}[0]{\mathfrak{N}}
\newcommand{\FO}[0]{\mathfrak{O}} \newcommand{\FP}[0]{\mathfrak{P}}
\newcommand{\FQ}[0]{\mathfrak{Q}} \newcommand{\FR}[0]{\mathfrak{R}}
\newcommand{\FS}[0]{\mathfrak{S}} \newcommand{\FT}[0]{\mathfrak{T}}
\newcommand{\FU}[0]{\mathfrak{U}} \newcommand{\FV}[0]{\mathfrak{V}}
\newcommand{\FW}[0]{\mathfrak{W}} \newcommand{\FX}[0]{\mathfrak{X}}
\newcommand{\FY}[0]{\mathfrak{Y}} \newcommand{\FZ}[0]{\mathfrak{Z}}

%Modifikation der Variablen
\renewcommand{\phi}[0]{\varphi}
\newcommand{\eps}[0]{\varepsilon}

%zusÔøΩtzliche Features
\newcommand{\quer}[0]{\overline}
\newcommand{\id}[0]{\operatorname{id}}		% IdentitÔøΩt
\newcommand{\Sp}[0]{\operatorname{Sp}}		% Spektrum eines Elements
\newcommand{\eins}[0]{\mathbf{1}}			% Eine Eins in allgemeinerem Kontext, z.B. in einem Ring
\newcommand{\auf}[1]{\quad\stackrel{#1}{\longrightarrow}\quad}
\newcommand{\ad}[0]{\operatorname{Ad}}
\newcommand{\ext}[0]{\operatorname{Ext}}
\newcommand{\fin}[0]{{\subset\!\!\!\subset}}
\newcommand{\Hom}[0]{\operatorname{Hom}}
\newcommand{\Aut}[0]{\operatorname{Aut}}
\newcommand{\dimnuc}[0]{\dim_{\mathrm{nuc}}}
\newcommand{\dimrok}[0]{\dim_{\mathrm{Rok}}}
\newcommand{\dimnuceins}[0]{\dimnuc^{\!+1}}
\newcommand{\dimrokeins}[0]{\dimrok^{\!+1}}
\newcommand{\im}[0]{\operatorname{im}}
\newcommand{\dst}[0]{\displaystyle}
\newcommand{\cstar}[0]{$\mathrm{C}^*$}
\newcommand{\Ost}[0]{\CO_\infty^{\mathrm{st}}}

% theorems
\newtheorem{satz}{Satz}[section]		% <--- optional, zÔøΩhlt so mit den Abschnitten
\newtheorem{cor}[satz]{Corollary}
\newtheorem{lemma}[satz]{Lemma}
\newtheorem{prop}[satz]{Proposition}
\newtheorem{theorem}[satz]{Theorem}

\theoremstyle{definition}
\newtheorem{defi}[satz]{Definition}
\newtheorem{nota}[satz]{Notation}
\newtheorem{rem}[satz]{Remark}
\newtheorem{example}[satz]{Example}
\newtheorem{question}[satz]{Question}

\newenvironment{bew}{\begin{proof}[Proof]}{\end{proof}}

\begin{abstract}
We investigate outer symmetries on unital Kirchberg algebras with respect to the Rokhlin property and finite Rokhlin dimension. In stark contrast to the restrictiveness of the Rokhlin property, every such action has Rokhlin dimension at most one. A consequence of these observations is a relationship between the nuclear dimension of an $\CO_\infty$-absorbing \cstar-algebra and its $\CO_2$-stabilization. We also give a more direct and alternative approach to this result. Several applications of this relationship are discussed to cover a fairly large class of $\CO_\infty$-absorbing \cstar-algebras that turn out to have finite nuclear dimension.
\end{abstract}

\maketitle

%\thispagestyle{empty}
%\newpage \tableofcontents
%\setcounter{page}{1}

\setcounter{section}{-1}

\section{Introduction}
\noindent
Studying group actions on \cstar-algebras has become a more and more prominent field of research during the last decade. While the general task of classification is far too complicated to expect satisfying results, there begins to develop a nice classification theory for certain finite group actions on Kirchberg algebras. Notable progress is made currently within this area by Phillips, leaning also on Köhler's research regarding an equivariant version of the UCT. 
Izumi's work in this area (see \cite{Izumi} and \cite{Izumi2}) has been seminal for the theory as a whole and for the presented results within this paper in particular.

So far, one of the driving forces behind many attempts to classify such actions has been the Rokhlin property. 
A high degree of rigidity of actions with the Rokhlin property is one of the main reasons why such attempts could succeed. On the other hand, this causes severe obstructions to the existence of such actions, see \cite{Izumi}.
To circumvent these obstructions, Hirsh\-berg, Winter and Zacharias invented a higher rank version of the Rokhlin property in \cite{HWZ}, namely finite Rokhlin dimension. In this sense, the classical Rokhlin property has to be understood as Rokhlin dimension zero.

While finite Rokhlin dimension comes with a higher amount of flexibility, it has been shown to behave well with underlying \cstar-algebras of finite nuclear dimension. In fact, the first part of this paper serves to give an interesting class of examples that showcases how much more flexibility is possible. Namely, we take a look at what happens if one makes the jump from Rokhlin dimension 0 to Rokhlin dimension 1 in the case of $\IZ_2$-actions on unital Kirchberg algebras. It turns out that in fact every outer $\IZ_2$-action on a unital Kirchberg algebra has Rokhlin dimension at most one. This relies heavily on Goldstein's and Izumi's remarkable main result of \cite{GoldIz}, which asserts that every such outer action absorbs a faithful quasi-free action on $\CO_\infty$. Thus we obtain our result by proving that a faithful, quasi-free action on $\CO_\infty$ has Rokhlin dimension one.

The second part of this paper focuses on consequences related to the nuclear dimension of $\CO_\infty$-absorbing \cstar-algebras. An easily deduced consequence of the first part is that there is a relationship between the nuclear dimension of an $\CO_\infty$-absorbing \cstar-algebra and the nuclear dimension of its $\CO_2$-stabilization. Said dimensional relationship is described by an inequality that relates $\dimnuc(A\otimes\CO_\infty)$ to $\dimnuc(A\otimes\CO_2)$ for any \cstar-algebra $A$. The main content of the second part is to improve this estimate with a more direct approach, which is independent of the first part of the paper.

The usefulness of such a relationship is apparent considering that $\CO_2$-absorption is a very strong property, in fact so rigid that Kirchberg was able to classify $\CO_2$-absorbing \cstar-algebras via their prime ideal spaces. Note also that our abstract approach, unlike similar results as in \cite{WiZa, Enders, RSS} or \cite{RST}, neither makes use of any model systems for the \cstar-algebras in question, nor requires any UCT assumption. As an immediate application, we can recover the result in \cite[Section 7]{MatuiSato}, i.e.~that Kirchberg algebras have nuclear dimension at most 3. As a second application, we consider the non-simple case of continuous field \cstar-algebras with Kirchberg fibres. It turns out that, as conjectured in \cite{TiWi}, the nuclear dimension depends rather on the fibres than the complexity of the underlying topological space. For example, if the prime ideal space of such a continuous field \cstar-algebra is finite dimensional, compact and metrizable, then the nuclear dimension is finite.

The authors would like to thank the referee for several helpful suggestions that ended up making this paper more reader-friendly and self-contained.

\section{Preliminaries}

\begin{nota} Unless specified otherwise, we will stick to the following notation throughout the paper.
\begin{itemize}
\item $A$ denotes a separable \cstar-algebra.
\item $\alpha, \beta$ or $\gamma$ denote group actions on a \cstar-algebra. Mostly, the acting group will be $\IZ_2=\IZ/2\IZ$.
\item If $M$ is some set and $F\subset M$ some finite subset, we write $F\fin M$.
\item For $\eps>0$ und $a,b$ in some \cstar-algebra, we write $a=_\eps b$ for $\|a-b\|\leq\eps$.
\item $\CZ$ denotes the Jiang-Su algebra.
\item While $\dimnuc(A)$ denotes, as usual, the nuclear dimension for a \cstar-algebra $A$, we also use the (non-standard, but very convenient) notation $\dimnuceins(A)=\dimnuc(A)+1$.
\item Analogously, the same goes for Rokhlin dimension $\dimrok$ and $\dimrokeins$.
\end{itemize}
\end{nota}

Let us recall the notion of Rokhlin dimension (see \cite{HWZ}) for finite group actions.

\begin{defi} \label{dimrok} 
Let $G$ be a finite group, $A$ a unital \cstar-algebra and \linebreak $\alpha: G\curvearrowright A$ an action via automorphisms. $\alpha$ is said to have Rokhlin dim\-ension $d$, written $\dimrok(\alpha)=d$, if $d$ is the smallest natural number with the following property:

For all $\eps>0$ and $F\fin A$, there exist positive contractions 
$(f_g^{(l)})_{g\in G}^{l=0,\dots,d}$ satisfying\vspace{-1mm}
\begin{itemize}
\item[(1)] $\displaystyle \eins_A=_\eps\sum_{l=0}^d \sum_{g\in G} f_g^{(l)}$.
\item[(2)] $\alpha_g(f_h^{(l)})=_\eps f_{gh}^{(l)}$\quad for all $l=0,\dots,d$ and $g,h\in G$.\vspace{1mm}
\item[(3)] $\|f_g^{(l)}f_h^{(l)}\|\leq\eps$\quad for all $l=0,\dots,d$ and $g\neq h$ in $G$.\vspace{1mm}
\item[(4)] $\|[f_g^{(l)},a]\|\leq\eps$\quad for all $l=0,\dots,d$, $g\in G$ and $a\in F$.\vspace{1mm}
\end{itemize}
If there is no such $d$, we write $\dimrok(\alpha)=\infty$.
\end{defi}

The usefulness of this notion is for instance illustrated in the following result from \cite{HWZ}:

{\theorem \label{HWZ} For a finite group $G$, a unital \cstar-algebra $A$ and an action $\alpha: G\curvearrowright A$, we have
\[\dimnuceins(A\rtimes_\alpha G)\leq\dimrokeins(\alpha)\cdot\dimnuceins(A).\] }

%%%%%%%%%%%%%%
\section{Outer $\IZ_2$-actions on Kirchberg algebras}
\noindent
In this section, we prove that outer $\IZ_2$-actions on unital Kirchberg algebras always have Rokhlin dimension at most 1. This will have an interesting application in the third section.

{
\lemma \label{Rokhlin elements} 
Let $A$ be a unital Kirchberg algebra and let $p\in A$ be a projection with vanishing $K_0$-class. Let $u=\eins-2p$. For all $\eps>0$, there exist positive contractions $f^{(0)}$ and $f^{(1)}$ in $A$ such that
\[ f^{(i)}\perp uf^{(i)}u^*\quad\text{for}\quad i=0,1\quad\text{and}\quad \sum_{i=0,1} f^{(i)}+uf^{(i)}u^* =_{\eps} \eins. \]
}
\begin{proof}
Let $n$ be an odd natural number with $n\geq\frac{1}{\eps}$. We can choose (using \cite{Cuntz}) projections $q_1,\dots,q_n$ and $r_1,\dots,r_n$ in $A$ such that
$$p=r_1+\dots+r_n,\quad \eins-p = q_1+\dots+q_n$$
and
$$[r_j]_0 = 0,\quad [q_j]_0=\begin{cases}1 &,\quad j~\text{odd}\\ -1 &,\quad j~\text{even} \end{cases}$$
for all $j=1,\dots,n$. It follows that
$$u(q_j+q_{j+1})=q_j+q_{j+1}\quad\text{and}\quad u(r_j+r_{j+1})=-(r_j+r_{j+1})$$
for all $j=1,\dots,n-1.$
Since $[q_j+q_{j+1}]_0=0=[r_j+r_{j+1}]_0$ for all $j=1,\dots,n-1$, it follows from \cite{Cuntz} that these projections are Murray-von-Neumann equivalent inside $A$. Hence, there is a partial isometry $v_j\in A$ such that $v_j^*v_j=q_j+q_{j+1}$ and $v_jv_j^*=r_j+r_{j+1}$. Setting
$$e_j = \frac{1}{2}( |v_j|+|v_j^*|+v_j+v_j^*)$$
yields projections with
$$e_j+ue_ju^* = p_j+p_{j+1}+r_j+r_{j+1}\quad\text{for all}~j=1,\dots,n-1.$$
Now we define
$$f^{(0)} = \sum_{j=2\atop j~\text{even}}^{n-1} \frac{j}{n} e_j\quad\text{and}\quad f^{(1)} = \sum_{j=1\atop j~\text{odd}}^{n-2} \frac{n-j}{n} e_j.$$
Because of the properties of the $e_j$, we observe $f^{(i)}\perp uf^{(i)}u^*$ for $i=0,1$. Moreover, one also gets
$$f^{(0)}+uf^{(0)}u^*+f^{(1)}+uf^{(1)}u^*=_{\frac{1}{n}} \eins.$$
This finishes the proof.
\end{proof}

\cor \label{quasifree} A faithful, quasi-free action $\IZ_2\curvearrowright\CO_\infty$ in the sense of \cite{GoldIz} has Rokhlin dimension 1.
\begin{proof}
By \cite[5.2]{GoldIz}, there is only one such action up to conjugacy. It follows from \cite[Section 6]{GoldIz}, and particularly from the proof of \cite[6.2]{GoldIz}, that for any non-trivial projection $p\in\CO_\infty$ with $[p]_0=0\in K_0(\CO_\infty)$, considering the unitary $u=\eins-2p$, the action
$$\gamma: \IZ_2\curvearrowright \CO_\infty\cong\bigotimes_\IN \CO_\infty,\quad \gamma=\bigotimes_\IN \ad(u)$$
is faithful and quasi-free. By \ref{Rokhlin elements}, it is clear that $\gamma$ has Rokhlin dimension at most 1. On the other hand, it is known that no finite group action on $\CO_\infty$ can have the Rokhlin property. Hence the proof is complete.
\end{proof}

\theorem \label{rokhlin} Let $A$ be a unital Kirchberg algebra and $\alpha: \IZ_2\curvearrowright A$ an action. If $\alpha$ is outer, then $\alpha$ has Rokhlin dimension at most 1.
\begin{proof}
Let $\gamma$ be a faithful, quasi-free action of $\IZ_2$ on $\CO_\infty$. By \cite[5.1]{GoldIz}, $\alpha$ is conjugate to $\alpha\otimes\gamma$. Hence the statement follows immediately from \ref{quasifree}.
\end{proof}

\defi \label{O standard} Let $p\in\CO_\infty$ be a projection with trivial class in $K_0(\CO_\infty)$. Define $\Ost := p\CO_\infty p$. 

\begin{rem} \label{beta} 
In \cite[4.7]{Izumi}, it was shown that there exists an action $\beta: \IZ_2\curvearrowright\CO_2$ such that $\CO_2\rtimes_\beta\IZ_2\cong\Ost\otimes M_{2^\infty}$.
\end{rem}
 
It follows from \ref{rokhlin} that this action has Rokhlin dimension 1. We will use this in the next section to deduce that for every unital \cstar-algebra $A$, if $A\otimes\CO_2$ has finite nuclear dimension, then so does $A\otimes\CO_\infty$.

%%%%%%%%
\section{The nuclear dimension of $\CO_\infty$-absorbing \cstar-algebras}
\noindent
With some additional observations, we use \ref{beta} to derive a dimensional inequality between an $\CO_\infty$-absorbing \cstar-algebra and its $\CO_2$-stabilization.

{
\lemma[see {\cite[Section 5]{MatuiSato}}] \label{ZUHF} 
For any \cstar-algebra $A$ and for any UHF algebra $\CU$ of infinite type, we have the inequality
\[
\dimnuceins(A\otimes\CZ)\leq 2\dimnuceins(A\otimes\CU).
\]
}
\begin{proof}
In \cite[Section 5]{MatuiSato}, the following fact was established: For every $\eps>0$, there exists $N\in\IN$, such that for all $m\geq N$, there exist two c.p.c.~order zero maps $\psi_0,\psi_1: M_m\to\CZ$ such that $\psi_0(\eins)+\psi_1(\eins)=_\eps\eins_\CZ$. With a standard argument, we may thus obtain two c.p.c.~order zero maps
\[
\phi_0,\phi_1: \CU\to\CZ_\infty\quad\text{with}\quad \phi_0(\eins_\CU)+\phi_1(\eins_\CU) = \eins_\CZ.
\]
Let $A$ be a \cstar-algebra and $d=\dimnuc(A\otimes\CU)$. It suffices to check that the nuclear dimension of the canonical embedding $\id_A\otimes\eins: A\to A\otimes\CZ$ is at most $2d+1$. We remark that the nuclear dimension of this map is the same as the nuclear dimension of its composition with the canonical embedding $A\otimes\CZ\to (A\otimes\CZ)_\infty$. This follows easily from the projectivity of finite dimensional \cstar-algebras with respect to c.p.c.~order zero maps, see \cite[3.10]{SWZ} for a very similar statement.

Now let $F\fin A$ be a finite subset and $\eps>0$. Find a finite dimensional \cstar-algebra $\CF$, a c.p.c.~map $\psi: A\to \CF$ and c.p.c.~order zero maps $\kappa_0,\dots,\kappa_d: \CF\to A\otimes\CU$ such that
\[
\sum_{j=0}^d \kappa_j\circ\psi(x)=_\eps x\otimes\eins_\CU\quad\text{for all}~x\in F.
\]
It thus follows that
\[
\sum_{j=0}^d\sum_{i=0,1} (\id_A\otimes\phi_i)\circ\kappa_j\circ\psi(x) =_\eps x\otimes\eins_\CZ\quad\text{for all}~x\in F.
\]
Thus we obtain a $2(d+1)$-decomposable approximation
\[
\xymatrix{
A \ar[rrrr]^{\id_A\otimes\eins_\CZ} \ar[rrd]_\psi &&&& (A\otimes\CZ)_\infty \\
&& \CF \ar[rru]_{\hspace{8mm}\sum_{j=0}^d\sum_{i=0,1} (\id_A\otimes\phi_i)\circ\kappa_j} &&
}
\]
for $F$ up to $\eps$. This finishes the proof.
\end{proof}

{\cor For any unital \cstar-algebra $A$, we have
\[\dimnuceins(A\otimes\CO_\infty)\leq 4\dimnuceins(A\otimes\CO_2).\]}
\begin{proof}
By \ref{beta}, there is an action $\beta: \IZ_2\curvearrowright\CO_2$ with Rokhlin dimension 1 such that $\CO_2\rtimes_\beta\IZ_2\cong\Ost\otimes M_{2^\infty}$. Define the action $\alpha=\id_A\otimes\beta: \IZ_2\curvearrowright A\otimes\CO_2$. Then
\[ A\otimes\Ost\otimes M_{2^\infty} \cong (A\otimes\CO_2)\rtimes_\alpha\IZ_2.\]
Moreover, it is clear that $\dimrok(\alpha)\leq\dimrok(\beta)=1$. Since $\Ost\otimes\CZ\cong\Ost$, we get
\[\begin{array}{lll}
\dimnuceins(A\otimes\Ost) &\stackrel{\ref{ZUHF}}{\leq}& 2\dimnuceins(A\otimes\Ost\otimes M_{2^\infty}) \vspace{3mm}\\
&=& 2\dimnuceins((A\otimes\CO_2)\rtimes_\alpha\IZ_2) \vspace{2mm}\\
&\stackrel{\ref{HWZ}}{\leq}& 2\dimrokeins(\alpha)\dimnuceins(A\otimes\CO_2)\vspace{3mm}\\
&\leq& 4\dimnuceins(A\otimes\CO_2).
\end{array}\]
By Brown's theorem, $\Ost$ is stably isomorphic to $\CO_\infty$, see \ref{O standard}. Since nuclear dimension is invariant under stable isomorphism (see \cite[2.8]{WiZa}), we are done.
\end{proof}

For the rest of this section, we would like to present a more direct approach to this dimensional inequality, not using group actions. This even allows us to improve the final estimate.

{
\theorem \label{dimension} For any \cstar-algebra $A$, we have the inequality
\[\dimnuceins(A\otimes\CO_{\infty})\leq 2 \dimnuceins(A\otimes \CO_{2}).\]
}
\begin{bew}
Since $\CO_{\infty}$ is strongly self-absorbing and $\CO_{2}$ absorbs $\CO_{\infty}$, we may assume $A\cong A\otimes\CO_{\infty}$. Let $F\fin A$ be a finite set of contractions, $\eps>0$ and let $h\in \CO_{\infty}$ be a positive element with spectrum equal to $[0,1]$. Since $\CO_{\infty}$ has real rank zero (see \cite[4.1.1]{Rordam2001}), we can find $k\in \IN$, mutually orthogonal non-trivial projections $p_{1},\ldots,p_{k}\in \CO_{\infty}$ and numbers $0\leq\lambda_{1}\leq\ldots\leq\lambda_{k}\leq1$ such that
\[
h =_\eps \sum_{i=1}^{k}\lambda_{i}p_{i}.
\]
Observe that we automatically have $\lambda_{1}\leq \eps$ and $\lambda_{i+1}-\lambda_{i}\leq\eps$ for $i=1,\ldots,k-1$.

Find a projection $0\neq q_{k}\leq p_{k}$ with $[q_{k}]_{0}=0\in K_{0}(\CO_{\infty})$ and set $\tilde{p}_{k}:=q_{k}$.
Let $0\leq j\leq k-2 $ and assume that we have already constructed $\tilde{p}_{k},\ldots,\tilde{p}_{k-j}$. Find a non-trivial projection $q_{k-j-1}\leq p_{k-j-1}$ with $[q_{k-j-1}]_{0}+[p_{k-j}-\tilde{p}_{k-j}]_{0}=0\in K_{0}(\CO_{\infty})$, and define
\[{\tilde{p}_{k-j-1}}:=p_{k-j}-\tilde{p}_{k-j}+q_{k-j-1}.\]
By construction, all $\tilde{p}_{i}$ have trivial class in $K_{0}(\CO_{\infty})$. If $\tilde{h}_{0}:=\sum_{i=1}^{k}\lambda_{i}\tilde{p}_{i}$, we can conclude
\[\begin{array}{lcl} \displaystyle\sum_{i=1}^{k}\lambda_{i}p_{i} - \tilde{h}_{0}&=& \displaystyle \sum_{i=1}^{k}\lambda_{i}p_{i}-\sum_{i=1}^{k}\lambda_{i}\tilde{p}_{i}\\\\
&=& \displaystyle\lambda_{1}(p_1-q_1)+\sum_{i=1}^{k-1}(\lambda_{i+1}-\lambda_{i})(p_{i+1}-q_{i+1}).
\end{array}\]
This implies $h=_{2\eps}\tilde{h}_{0}$
and gives an approximation of $h$ by a positive linear combination of pairwise orthogonal projections with trivial $K_{0}$-classes. By \cite[4.2.3]{Rordam2001}, $\CO_{2}$ embeds unitally into the corner $\tilde{p}_{i}\CO_{\infty}\tilde{p}_{i}$ for each $i$. As $\CO_{\infty}$ is purely infinite and simple, all $\tilde{p}_{i}$ are Murray-von Neumann equivalent, see \cite[1.4]{Cuntz}. With $k_0=k$, the above embeddings extend to an injective $*$-homo\-morphism
\[
\tilde{\iota}_{0} : M_{k_{0}}\otimes\CO_{2}\rightarrow \CO_{\infty}\quad\text{with}\quad \tilde{\iota}_{0}(e_{ii}\otimes\eins)=\tilde{p}_i.
\]
Clearly, $\tilde{h}_{0}$ is contained in the image of $\tilde{\iota}_{0}$. Analogous considerations for $\eins-h$ lead to a natural number $k_1$, an injective $*$-homomorphism
$\tilde{\iota}_{1}:M_{k_{1}}\otimes\CO_{2}\rightarrow \CO_{\infty}$
and a positive contraction $\tilde{h}_{1}\in \im(\tilde{\iota}_{1})$ satisfying $\eins-h=_{2\eps} \tilde{h}_{1}$.

Since the statement is trivial in the case that the nuclear dimension of $A\otimes \CO_{2}$ is infinite, we may assume $\dimnuc(A\otimes \CO_{2})=d<\infty$. Since $\CO_{2}\cong M_{k_{i}}\otimes \CO_{2}$ for $i=0,1$, let $\iota_{i}:\CO_{2}\rightarrow \CO_{\infty}$ be the injective $*$-homomorphism corresponding to $\tilde{\iota}_{i}$ and set $h_{i}:=\iota_{i}^{-1}(\tilde{h}_{i})\in\CO_{2}$ for $i=0,1$.

There exist finite dimensional \cstar-algebras $E_{i,0},\ldots, E_{i,d}$, c.p.c.~maps $\eta_{i}:A\otimes\CO_{2}\rightarrow E_{i,0}\oplus\ldots\oplus E_{i,d}$ and c.p.c.~order zero maps $\kappa_{i,j}:E_{i,j}\rightarrow A\otimes\CO_{2}$ such that for all $a\in F$
\[(\star)\qquad
\sum_{j=0}^{d}\kappa_{i,j}\circ \eta_{i}(a\otimes h_{i}) =_\eps a\otimes h_{i}\]
holds. Define a c.p.c.~map
\[\psi : A\rightarrow \bigoplus_{i=0,1}\bigoplus^{d}_{j=0} E_{i,j}\quad\text{via}\quad
 \psi =(\eta_{0}\oplus \eta_{1})\circ (\id_{A}\otimes h_{0}\oplus \id_{A}\otimes h_{1}),\]
 and c.p.c.~order zero maps 
$\tilde{\phi}_{i,j} : E_{i,j}\rightarrow A\otimes \CO_{\infty}$ via
$\tilde{\phi}_{i,j} :=(\id_{A}\otimes\iota_{i})\circ\kappa_{i,j}.$
Using the estimate $(\star)$, we get
\[
\begin{array}{lll}
\dst \Bigl(\sum_{i=0,1}\sum_{j=0}^{d} \tilde{\phi}_{i,j}\Bigl)\circ\psi(a) &=& \dst \sum_{i=0,1}\id_{A}\otimes \iota_{i}\Bigl(\sum_{j=0}^{d} \kappa_{i,j}\circ\eta_i(a\otimes h_{i})\Bigl)\vspace{2mm}\\
&=_{2\eps}&\dst \sum_{i=0,1}\id_{A}\otimes\iota_{i}(a\otimes h_{i}) \vspace{2mm}\\
&=& a\otimes(\tilde{h}_{0}+\tilde{h}_{1}) \vspace{2mm}\\
&=_{4\eps}& a\otimes\eins.
\end{array}\]
for all $a\in F$. As $A$ is $\CO_{\infty}$-absorbing, \cite[2.7]{TomsWinter2007} implies that there exists a $*$-homomorphism $\sigma:A\otimes \CO_{\infty}\rightarrow A$ satisfying
\[\sigma(a\otimes \eins)=_\eps a\quad\text{for all}\; a\in F.\]
Define c.p.c.~order zero maps $\phi_{i,j}:=\sigma\circ\tilde{\phi}_{i,j}$. We get
\[
\Bigl(\sum_{i=0,1} \sum_{j=0}^{d} \phi_{i,j}\Bigl)\circ\psi(a) = \sigma\left[\Bigl(\sum_{i=0,1}\sum_{j=0}^{d} \tilde{\phi}_{i,j}\Bigl)\circ\psi(a)\right]=_{6\eps} \sigma(a\otimes\eins) =_\eps a
\]
for all $a\in F$. Since $F\fin A$ and $\eps>0$ were arbitrary, this shows that $\dimnuc(A\otimes\CO_{\infty})\leq 2(d+1)-1$ and completes the proof.
\end{bew}

This reproduces the result of \cite[Section 7]{MatuiSato} about finite nuclear dimension of Kirchberg algebras. Note that no UCT assumption is needed.

{\cor For every Kirchberg algebra $A$, we have $\dimnuc(A)\leq 3$.}
\begin{proof}
Because of \cite[2.8]{WiZa}, we may assume that $A$ is unital. It has also been shown in \cite{WiZa} that $\dimnuc(\CO_2)=1$. The assertion now follows directly from \ref{dimension} and Kirchberg's absorption theorems (see \cite{KiPhi})
\[A\cong A\otimes\CO_\infty\quad\text{and}\quad A\otimes\CO_2\cong\CO_2.\]
\end{proof}

One of the strengths of \ref{dimension} lies in the fact that no simplicity assumptions are needed. In fact, there are plenty of highly non-simple examples that we can treat. This is because $\CO_2$-absorbing \cstar-algebras make up a very small class compared to $\CO_\infty$-absorbing \cstar-algebras. We can illustrate this nicely on bundles.

{\theorem \label{bundles} Let $X$ be a compact metrizable space. Let $A$ be an $\CO_\infty$-absorbing, continuous $\CC(X)$-algebra whose fibres are Kirchberg algebras. Then $A$ has finite nuclear dimension. (In fact, we have that $\dimnuc(A)\leq 7.$)
}
\begin{proof}
If we apply Kirchberg's classification of $\CO_2$-absorbing \cstar-algebras by their prime ideal spaces, we get $A\otimes\CO_2 \cong \CC(X)\otimes\CO_2.$

Now it follows from \cite[3.7]{KiRo} that there is a one-dimensional, compact metric space $Y$ and an embedding $\iota: \CC(Y)\to\CC(X)\otimes\CO_2$ such that 
\[
\CC(X)\otimes\eins_{\CO_2}\subset\iota(\CC(Y))\subset\CC(X)\otimes\CO_2.
\]
In particular, the nuclear dimension of the canonical embedding $\id_{\CC(X)}\otimes\eins_{\CO_2}: \CC(X)\to\CC(X)\otimes\CO_2$ is at most one. From this, it follows that
\[
\begin{array}{ccl}
\dimnuceins(\CC(X)\otimes\CO_2) &=& \dimnuceins(\id_{\CC(X)\otimes\CO_2}\otimes\eins_{\CO_2}) \\\\
& \leq & \dimnuceins(\id_{\CC(X)}\otimes\eins_{\CO_2})\cdot\dimnuceins(\CO_2) \\\\
& = & 2\cdot 2 = 4.
\end{array}
\]
In particular, this yields
\[\dimnuceins(A)=\dimnuceins(A\otimes\CO_\infty)\stackrel{\ref{dimension}}{\leq} 2\dimnuceins(A\otimes\CO_2) = 8.\]
\end{proof}

\cor Let $X$ be a compact metrizable space of finite covering dimension. Let $A$ be a continuous $\CC(X)$-algebra whose fibres are Kirchberg algebras. Then $A$ has finite nuclear dimension.
\begin{proof}
It follows from \cite[5.11]{BlaKi} that $A$ is automatically $\CO_\infty$-stable. Apply \ref{bundles} to get the statement.
\end{proof}

\question Do nuclear, $\CO_2$-absorbing \cstar-algebras always have finite nuclear dimension? In particular, the case of non-Hausdorff primitive ideal spaces is an intriguing open problem.

\end{document}